# EDGE-PROMOTING RECONSTRUCTION OF ABSORPTION AND DIFFUSIVITY IN OPTICAL TOMOGRAPHY

A. HANNUKAINEN[†], L. HARHANEN[‡], N. HYVÖNEN[†], AND H. MAJANDER[§]

**Abstract.** In optical tomography a physical body is illuminated with near-infrared light and the resulting outward photon flux is measured at the object boundary. The goal is to reconstruct internal optical properties of the body, such as absorption and diffusivity. In this work, it is assumed that the imaged object is composed of an approximately homogeneous background with clearly distinguishable embedded inhomogeneities. An algorithm for finding the *maximum a posteriori estimate* for the absorption and diffusion coefficients is introduced assuming an edge-preferring prior and an additive Gaussian measurement noise model. The method is based on iteratively combining a lagged diffusivity step and a linearization of the measurement model of diffuse optical tomography with priorconditioned LSQR. The performance of the reconstruction technique is tested via three-dimensional numerical experiments with simulated measurement data.

**1. Introduction.** The objective of *optical tomography* (OT) is to deduce useful information about the internal structure of a physical body by illuminating it with near-infrared light and measuring the corresponding outward photon flux at a finite number of sensors on the surface of the body. To be more precise, the standard aim is to reconstruct the absorption and (reduced) scattering coefficients inside the imaged object. The potential applications of OT include screening for breast cancer and mapping the structure and function in a neonatal brain. For medical and instrumental details of OT, we refer to the review articles [2, 3, 5, 12, 17, 22] and the references therein.

We model the light propagation in (strongly scattering) tissue by the *diffusion approximation* (DA) of the *radiative transfer equation* (RTE), which can be considered one of the standard choices in medical imaging [2]. See, e.g., [30] for employing the RTE directly as the forward model for OT. In the frequency domain, the forward problem associated to the DA consists of an (elliptic, steady state) diffusion equation accompanied by a Robin boundary condition that corresponds to the inward photon flux at the object boundary. The two coefficients in the diffusion equation, i.e. the diffusivity and the (complexified) absorption, can be represented as simple functions of the absorption, the scatter and the harmonic modulation frequency of the input boundary flux. For simplicity, we choose the diffusion and absorption coefficients as the to-be-reconstructed parameters, but it would be straightforward to adapt our algorithm so that the scatter was directly included as one of the two primary unknowns.

In this work, we assume it is *a priori* known that the material parameters of the imaged object are approximately constant apart from a few clearly distinguishable inclusions. Such a setting could be encountered, e.g., when detecting cancerous anomalies within healthy tissue. We quantify this auxiliary information by introducing Perona–Malik [25] prior probability densities for the logarithms of the discretized diffusion and absorption coefficients. Under the additional assumption of additive Gaussian measurement noise, we deduce the posterior density, i.e., the joint conditional probability density of the parameters given the noisy measurements. Then the computation of the *maximum a posteriori* (MAP) estimate for the unknown coefficients corresponds to finding the minimizer of a certain Tikhonov functional that is nonquadratic in both the discrepancy and the penalty term. The method for reconstructing the internal conductivity of a physical body by *electrical impedance tomography* (EIT) considered in [15] is based on iteratively minimizing such a Tikhonov functional; the main goal of this work is to generalize the ideas of [15] in order to introduce an algorithm that is capable of simultaneously reconstructing the diffusion and absorption coefficients in realistic three-dimensional OT.

[†]Aalto University, Department of Mathematics and Systems Analysis, P.O. Box 11100, FI-00076 Aalto, Finland (antti.hannukainen@aalto.fi, nuutti.hyvonen@aalto.fi).
[‡]Technical University of Denmark, Department of Applied Mathematics and Computer Science, Asmussens Alle, Building 322, DK-2800 Kgs. Lyngby, Denmark (lhar@dtu.dk). Aalto University, Department of Engineering Design and Production, P.O. Box 11000, FI-00076 Aalto, Finland.
[§]Aalto University, Department of Mathematics and Systems Analysis, P.O. Box 11100, FI-00076 Aalto, Finland (helle.majander@aalto.fi). École Polytechnique, Centre de Mathématiques Appliquées, Route de Saclay, 91128 Palaiseau Cedex France.

Our algorithm is composed of two nested iterations. The outer loop corresponds to introducing a standard quadratic Tikhonov functional that (heuristically) approximates the original nonquadratic one around the current iterate. The approximation is based on linearizing the measurement map corresponding to the DA and applying a single lagged diffusivity step [31] to get rid of the nonlinearity in the penalty term. The interior loop tackles a given quadratic Tikhonov functional by resorting to LSQR [23, 24] with priorconditioning: The positive definite matrix defining the penalty term is interpreted as the inverse covariance matrix of a 'linearized prior' in the spirit of priorconditioning [6, 7, 8, 9], and it is then employed as a (symmetric) preconditioner for the corresponding normal equation. Subsequently, the penalty term is dropped and the LSQR algorithm from [1] is applied to the resulting reduced normal equation. The leading idea of the preconditioning trick is to implicitly include the prior information about the inclusion boundaries from the previous iterate of the outer loop in the Krylov subspace structure of the LSQR phase. (Notice that priorconditioning is related to the transformation of a quadratic Tikhonov functional into the standard form [11, 14, 19].) See [28] for a related work where the linearized problem has been directly tackled by a Krylov method.

If it is *a priori* known that the inclusions only affect one of the two coefficients in the diffusion equation, the proposed algorithm produces good quality reconstructions from simulated noisy data independently of whether the inward photon flux at the object boundary is harmonically modulated or not. However, if there are anomalies in both the diffusion and the absorption coefficient, the performance of the algorithm is considerably better when the input is time-harmonic: For static inward boundary flux, the reconstructions exhibit a severe 'cross-talk' between the two coefficients. This observation, made also in the context of time-resolved measurements [29], is in line with the conclusion of [4]; see [16] for a twist, though. A special feature of our method is that it produces three-dimensional reconstructions with about $10^5$ degrees of freedom in only about ten minutes on a standard laptop. In addition, the algorithm does not require the knowledge of the background parameter levels.

This manuscript is organized as follows. Section 2 recalls the DA of the RTE and discusses an efficient way of sampling the derivatives of the associated measurement map with respect to the diffusion and absorption coefficients. Our preferred (finite-dimensional) Bayesian framework as well as the prior probability densities for the coefficients of the diffusion equation are introduced in Section 3. The actual reconstruction algorithm and its implementation are considered in Section 4. Finally, Section 5 presents the three-dimensional numerical examples.

**2. Measurement model.** In this work, we consider the forward model corresponding to the DA of the RTE. We assume the measurements are modulated with a fixed harmonic frequency $\omega \in \mathbb{R}$ and model the examined physical body as a bounded domain $\Omega \subset \mathbb{R}^3$ with a connected complement and a Lipschitz boundary. The photon sources on the object boundary are represented by the connected, mutually disjoint, open subsets $s_k \subset \partial\Omega$, $k = 1, \ldots, K$. The amplitude of the input flux $\Phi_k : \partial\Omega \to \mathbb{R}$ through each source is modelled by the corresponding characteristic function, i.e.,

$$\Phi_k(x) = \begin{cases} 1 & \text{if } x \in s_k, \\ 0 & \text{otherwise.} \end{cases} \tag{2.1}$$

We assume $\Omega$ is isotropic and denote the associated diffusion and absorption coefficients by $\kappa, \mu \in L_+^\infty(\Omega, \mathbb{R})$, respectively, with the physically reasonable definition

$$L_+^\infty(\Omega, \mathbb{R}) = \{v \in L^\infty(\Omega, \mathbb{R}) \,|\, \operatorname{ess\,inf} v > 0\}.$$

Take note that apart from $L^\infty(\Omega, \mathbb{R})$ and $L_+^\infty(\Omega, \mathbb{R})$ all other function spaces considered in this work have $\mathbb{C}$ as the associated multiplier field.

According to the DA [2, 18], the photon density $\varphi_k \in H^1(\Omega)$ corresponding to photon flux through the $k$th source is the unique solution [13] of the elliptic boundary value problem

$$\begin{cases} -\nabla \cdot (\kappa \nabla \varphi_k) + \left(\mu + \dfrac{\omega}{c}\mathrm{i}\right)\varphi_k = 0 & \text{in } \Omega, \\ \gamma\varphi_k + \dfrac{1}{2}\nu \cdot \kappa \nabla \varphi_k = \Phi_k & \text{on } \partial\Omega, \end{cases} \tag{2.2}$$

where $\nu : \partial\Omega \to \mathbb{R}^n$ is the exterior unit normal of $\partial\Omega$, $c$ is the speed of light that is assumed to be constant in $\Omega$, and the dimension-dependent constant $\gamma$ takes the value $\gamma = 1/4$ in our three-dimensional setting. Notice that if $\partial\Omega$ is of the Hölder class $\mathscr{C}^{1,1}$ and the boundaries of the sources $\partial s_k$, $k = 1, \ldots, K$, are of the class $\mathscr{C}^1$, then for any $\epsilon > 0$ the photon flux $\Phi_k$ belongs to $H^{1/2-\epsilon}(\partial\Omega)$ by the properties of zero continuation for Sobolev spaces (cf., e.g., [32]) and $\varphi_k \in H^{2-\epsilon}(\Omega)$ due to the standard regularity theory for elliptic boundary value problems [13].

Corresponding to the input flux through each source, we measure the net flux coming out of the object at the (open and connected) measurement sensors $m_j \subset \partial\Omega$, $j = 1, \ldots, J$. It is assumed that the sensors do not overlap with the sources nor with each other. We model the 'device function' of the $j$th sensor by its characteristic function

$$\Psi_j(x) = \begin{cases} 1 & \text{if } x \in m_j, \\ 0 & \text{otherwise.} \end{cases} \tag{2.3}$$

With this convention, the boundary measurement at the sensor $m_j$ corresponding to input flux through the source $s_k$ can be written as (cf. [18])

$$M_{jk} = \int_{\partial\Omega} \Psi_j \left(\gamma\varphi_k - \frac{1}{2}\nu \cdot \kappa \nabla \varphi_k\right) dS = 2\gamma \int_{\partial\Omega} \Psi_j \varphi_k \, dS. \tag{2.4}$$

Here the latter equality follows from the boundary condition in (2.2) together with the assumption that the supports of $\Psi_j$ and $\Phi_k$ are disjoint. For simplicity, we have ignored reflections at $s_k$ and $m_j$ in (2.2) and (2.4), respectively. Moreover, there certainly exist more case-specific models for the input amplitudes $\Phi_k$ and the device functions $\Psi_j$, but we have decided to model them as characteristic functions as this is arguably the most generic choice. Be that as it may, we include $\Phi_k$ and $\Psi_j$ explicitly in the following analysis in order to demonstrate that one could as easily resort to more complicated models.

The inverse problem we are interested in is to reconstruct useful information about both the diffusivity $\kappa$ and the absorption $\mu$ inside $\Omega$ based on noisy versions of the measurements $M_{jk}(\kappa,\mu) \in \mathbb{C}$, where $j = 1, \ldots, J$ and $k = 1, \ldots, K$.

We complete this section by recalling (cf., e.g., [2, Section 5.2]) an efficient way of sampling the Fréchet derivative of the map $(\kappa,\mu) \mapsto M_{jk}(\kappa,\mu)$. To this end, we need to introduce the 'dual problem' of (2.2), that is,

$$\begin{cases} -\nabla \cdot (\kappa \nabla \psi_j) + \left(\mu + \frac{\omega}{c}\mathrm{i}\right) \psi_j = 0 & \text{in } \Omega, \\ \gamma \psi_j + \frac{1}{2}\nu \cdot \kappa \nabla \psi_j = \Psi_j & \text{on } \partial\Omega, \end{cases} \qquad (2.5)$$

which has a unique solution $\psi_j \in H^1(\Omega)$ [13].

THEOREM 2.1. *The Fréchet derivative of the map*

$$\left[L_+^\infty(\Omega)\right]^2 \ni (\kappa,\mu) \mapsto M_{jk}(\kappa,\mu) \in \mathbb{C} \qquad (2.6)$$

*at $(\kappa,\mu) \in [L_+^\infty(\Omega)]^2$ in the direction $(\vartheta,\theta) \in [L^\infty(\Omega)]^2$ can be evaluated via*

$$\left(M'_{jk}(\kappa,\mu)\right)(\vartheta,\theta) = -\gamma \int_\Omega \vartheta \, \nabla \psi_j \cdot \nabla \varphi_k \, dx - \gamma \int_\Omega \theta \, \psi_j \varphi_k \, dx, \qquad (2.7)$$

*where $\varphi_k, \psi_j \in H^1(\Omega)$ are the solutions of (2.2) and (2.5), respectively.*

*Proof.* It is well known that the variational formulation of (2.5) is to find $\psi_j \in H^1(\Omega)$ such that (cf., e.g., [2, 13])

$$\int_\Omega \left(\kappa \nabla \psi_j \cdot \nabla v + \left(\mu + \frac{\omega}{c}\mathrm{i}\right) \psi_j v\right) dx + 2\gamma \int_{\partial\Omega} \psi_j v \, dS \;=\; 2 \int_{\partial\Omega} \Psi_j v \, dS \qquad (2.8)$$

for all $v \in H^1(\Omega)$.[1] It is also common knowledge that the map

$$\left[L_+^\infty(\Omega)\right]^2 \ni (\kappa,\mu) \mapsto \varphi_k \in H^1(\Omega),$$

where $\varphi_k = \varphi_k(\kappa,\mu)$ is the solution of (2.2), is Fréchet differentiable. The corresponding derivative at $(\kappa,\mu) \in [L_+^\infty(\Omega)]^2$ in the direction $(\vartheta,\theta) \in [L^\infty(\Omega)]^2$ is given as the unique solution of the following 'Born approximation -like' variational problem (cf., e.g., [10]): Find $\varphi'_k = (\varphi'_k(\kappa,\mu))(\vartheta,\theta) \in H^1(\Omega)$ satisfying

$$\int_\Omega \left(\kappa \nabla \varphi'_k \cdot \nabla v + \left(\mu + \frac{\omega}{c}\mathrm{i}\right) \varphi'_k v\right) dx + 2\gamma \int_{\partial\Omega} \varphi'_k v \, dS$$
$$= - \int_\Omega \vartheta \, \nabla \varphi_k \cdot \nabla v \, dS - \int_\Omega \theta \, \varphi_k v \, dS \qquad (2.9)$$

for all $v \in H^1(\Omega)$.

---
[1] We have left out the 'standard complex conjugation' of the test function as this makes the argumentation more straightforward.

Due to the trace theorem [13] and since the right-hand side of (2.4) is linear with respect to $\varphi_k$, it is straightforward to deduce that the Fréchet derivative of (2.6) can be evaluated through

$$\big(M'_{jk}(\kappa,\mu)\big)(\vartheta,\theta) \;=\; 2\gamma \int_\Omega \Psi_j \varphi'_k \, dS,$$

where $\varphi'_k \in H^1(\Omega)$ is the solution of (2.9). In consequence, choosing $v = \varphi'_k$ in (2.8), we may write

$$\frac{1}{\gamma} \big(M'_{jk}(\kappa,\mu)\big)(\vartheta,\theta) = \int_\Omega \Big(\kappa \nabla \psi_j \cdot \nabla \varphi'_k + \Big(\mu + \frac{\omega}{c} \mathrm{i}\Big)\psi_j \varphi'_k\Big) dx + 2\gamma \int_{\partial\Omega} \psi_j \varphi'_k \, dS$$

$$= - \int_\Omega \vartheta \, \nabla \psi_j \cdot \nabla \varphi_k \, dS - \int_\Omega \theta \, \psi_j \varphi_k \, dS,$$

where the second step follows from (2.9) with $v = \psi_j$. This completes the proof. □

In particular, observe that we do no need to solve (2.9) at any stage when computing derivatives of the measurements with respect to the to-be-reconstructed parameters $\kappa$ and $\mu$, but it is enough to substitute the (approximate) solutions of (2.2) and (2.5) together with the considered perturbation directions in (2.7).

**3. Bayesian framework and the choice of prior.** We now consider the discretized version of (2.2), modelling the diffusivity $\kappa$ and the absorption $\mu$ as exponential quantities,

$$\kappa(\varsigma) = \kappa_0 \exp(\varsigma) \qquad \text{and} \qquad \mu(\upsilon) = \mu_0 \exp(\upsilon), \tag{3.1}$$

where

$$\varsigma = \sum_{n=1}^N \varsigma_n \phi_n \qquad \text{and} \qquad \upsilon = \sum_{n=1}^N \upsilon_n \phi_n. \tag{3.2}$$

Here $\varsigma_n, \upsilon_n \in \mathbb{R}$, $n = 1, \ldots, N$, are real coefficients and $\phi_n \in H^1(\Omega)$, $n = 1, \ldots, N$, are the piecewise linear finite element basis functions corresponding to a chosen meshing of $\Omega$, excluding the nodes on a nonempty subset $S$ of $\partial\Omega$, which corresponds to assuming that $\kappa = \kappa_0$ and $\mu = \mu_0$ at those nodes. We denote by $\varsigma$ and $\upsilon$ both the corresponding vectors of coefficients, $\varsigma, \upsilon \in \mathbb{R}_+^N$ and the functions defined in (3.2); the exact meaning should be clear from the context. The positive real numbers $\kappa_0, \mu_0 > 0$ in (3.1) are the constant diffusivity and absorption levels that are the most compatible with the measurements; the determination of these constants as well as the selection of $S \subset \partial\Omega$ is explained in Section 4.3. According to our experience, reconstructing $\varsigma$ and $\upsilon$ produces significantly better results than directly considering the actual physical quantities $\kappa$ and $\mu$; this observation is in line with, e.g., the material in [28, 30]. The transformation also ensures the positivity of $\kappa$ and $\mu$.

We vectorize the measurements as functions of $\varsigma$ and $\upsilon$, i.e., set

$$M(\varsigma, \upsilon) = \big[M_{jk}\big(\kappa(\varsigma), \mu(\upsilon)\big)\big] \in \mathbb{C}^{JK},$$

and introduce the real-valued measurement vector

$$\mathcal{M}(\varsigma, \upsilon) = \begin{bmatrix} \mathrm{Re}\, M(\varsigma, \upsilon) \\ \mathrm{Im}\, M(\varsigma, \upsilon) \end{bmatrix} \in \mathbb{R}^{2JK}. \tag{3.3}$$

For $\mathcal{M}$, we assume the additive Gaussian noise model

$$\mathcal{V} = \mathcal{M}(\varsigma, \upsilon) + \eta, \qquad (3.4)$$

where $\eta \in \mathbb{R}^{2JK}$ is a realization of a normally distributed random variable with zero mean and a known, symmetric and positive definite covariance matrix $\Gamma \in \mathbb{R}^{2JK \times 2JK}$. It easily follows that the likelihood function, i.e. the probability density of the measurements given the parameters, is the same Gaussian but with a shifted mean,

$$p(\mathcal{V} \,|\, \varsigma, \upsilon) \propto \exp\Big(-\frac{1}{2}\big(\mathcal{V} - \mathcal{M}(\varsigma, \upsilon)\big)^{\mathrm{T}} \Gamma^{-1} \big(\mathcal{V} - \mathcal{M}(\varsigma, \upsilon)\big)\Big),$$

where the constant of proportionality is independent of $\varsigma$ and $\upsilon$.

We assume that the (discretized) logarithms of the diffusivity and absorption are *a priori* distributed according to the densities

$$p(\varsigma) \propto \exp\big(-a\,R(\varsigma)\big) \qquad \text{and} \qquad p(\upsilon) \propto \exp\big(-b\,R(\upsilon)\big) \qquad (3.5)$$

where $a, b > 0$ are free parameters and $R$ is of the form

$$R(u) = \int_\Omega r\big(|\nabla u(x)|\big)\, dx, \qquad (3.6)$$

with $r : \Omega \to \mathbb{R}_+$ being a suitable, continuously differentiable, monotonically increasing function that prefers edges over slow changes in the function $u$. In this work we resort exclusively to

$$r(t) = \frac{1}{2} T^2 \log\big(1 + (t/T)^2\big), \qquad (3.7)$$

which corresponds to Perona–Malik prior/regularization [25], with $T > 0$ being a small parameter that gives a threshold for detectable edges. Employing the standard TV prior [26] would lead to slightly more blurred reconstructions, whereas using the so-called TV$^q$ prior (cf., e.g., [20]) with, say, $q = 1/2$ would produce results that are qualitatively the same as the ones presented in Section 5. Since (3.5) gives the probability densities for the logarithms $\varsigma$ and $\upsilon$ of the physical quantities interest, the prior imposed on $\kappa$ and $\mu$ could be dubbed 'log-Perona–Malik', imitating the relation between normal and log-normal distributions.

Assuming $\varsigma$ and $\upsilon$ are independent, the Bayes' formula gives the posterior density

$$p(\varsigma, \upsilon \,|\, \mathcal{V}) \propto p(\mathcal{V} \,|\, \varsigma, \upsilon)\, p(\varsigma) p(\upsilon)$$

$$\propto \exp\Big(-\frac{1}{2}\big(\mathcal{V} - \mathcal{M}(\varsigma, \upsilon)\big)^{\mathrm{T}} \Gamma^{-1} \big(\mathcal{V} - \mathcal{M}(\varsigma, \upsilon)\big) - a\,R(\varsigma) - b\,R(\upsilon)\Big),$$

where the constants of proportionality do not depend on $\varsigma$ and $\upsilon$. Our aim now is to find the MAP estimate corresponding to this posterior, which is equivalent to finding the minimizer of the Tikhonov functional

$$F(\varsigma, \upsilon) := \frac{1}{2}\big(\mathcal{V} - \mathcal{M}(\varsigma, \upsilon)\big)^{\mathrm{T}} \Gamma^{-1} \big(\mathcal{V} - \mathcal{M}(\varsigma, \upsilon)\big) + a\,R(\varsigma) + b\,R(\upsilon). \qquad (3.8)$$

In what follows, we denote $\beta = [\varsigma^{\mathrm{T}}, \upsilon^{\mathrm{T}}]^{\mathrm{T}} \in \mathbb{R}^{2N}$ and, instead of $\mathcal{M}(\varsigma, \upsilon)$ and $F(\varsigma, \upsilon)$, we usually write $\mathcal{M}(\beta)$ and $F(\beta)$ for short.

**4. The algorithm.** The iterative method for minimizing (3.8) starts from the initial guess $\beta^{(0)} = 0 \in \mathbb{R}^{2N}$, which corresponds to the constant values $\kappa_0$ and $\mu_0$ for the diffusivity and the absorption, respectively (cf. (3.1)). We linearize $\mathcal{M}(\beta)$ around the current iterate $\beta^{(l)}$, write the necessary condition for minimizing the resulting Tikhonov functional, and take one *lagged diffusivity* step [31]. Altogether, this corresponds to finding the minimizer of a certain standard quadratic Tikhonov functional. The minimization problem is tackled by means of LSQR-aided priorconditioning, and the resulting parameter vector is taken to be the next iterate $\beta^{(l+1)}$.

In the following, we first briefly deduce the linear problem corresponding to the lagged diffusivity step and subsequently explain how we employ priorconditioning to approximately solve it. Finally, we summarize the whole reconstruction method as a complete algorithm in Section 4.3. For more information on the basic building blocks of the algorithm we refer to [15], where its variant was applied to EIT. Notice, however, that [15] considered the reconstruction of only one coefficient in an elliptic PDE, namely the conductivity, whereas here we are interested in both the diffusivity and the absorption. In particular, [15] did not employ an exponential model of the form (3.1) for the conductivity.

**4.1. Linearization and lagged diffusivity step.** Linearizing $\mathcal{M}(\beta)$ around $\beta^{(l)}$ in (3.8) leads to a new Tikhonov functional with a quadratic residual but non-quadratic penalty terms,

$$F^{(l)}(\beta) := \frac{1}{2}\big(y^{(l)} - J^{(l)}\beta\big)^{\mathrm{T}}\Gamma^{-1}\big(y^{(l)} - J^{(l)}\beta\big) + a\,R(\varsigma) + b\,R(\upsilon), \qquad (4.1)$$

where $\beta = [\varsigma^{\mathrm{T}}, \upsilon^{\mathrm{T}}]^{\mathrm{T}}$, the matrix $J^{(l)} \in \mathbb{R}^{2JK \times 2N}$ is the Jacobian of the map $\beta \mapsto \mathcal{M}(\beta)$ evaluated at $\beta^{(l)}$ and

$$y^{(l)} = \mathcal{V} - \mathcal{M}(\beta^{(l)}) + J^{(l)}\beta^{(l)} \in \mathbb{R}^{2JK}.$$

Let us briefly explain how to numerically compute $\mathcal{M}(\beta^{(l)})$ and $J^{(l)}$ needed in the above formulae: $\mathcal{M}(\beta^{(l)})$ is approximated by solving (2.2) with the coefficients (3.1) for the input fluxes $\Phi_k$, $k = 1, \ldots, K$, by a FEM, substituting the obtained solutions in (2.4) for all $j = 1, \ldots, J$ and separating the real and imaginary parts. On the other hand, the approximate elements of $J^{(l)}$ are obtained by plugging the aforementioned FEM solutions of (2.2) together with those of (2.5) in (2.7) and letting the perturbations $\vartheta$ and $\theta$ run in turns through the piecewise linear basis functions $\phi_n$, $n = 1, \ldots, N$, appearing in the representations (3.2). To be quite precise, this process only gives the elements in the Jacobian of the complex-valued measurements $M(\beta) \in \mathbb{C}^{JK}$ at $\beta^{(l)}$ with respect to the vector $[\exp(\beta^{(l)})] \in \mathbb{R}^{2N}$ (cf. (3.1)), but $J^{(l)}$ is subsequently obtained with the help of the chain rule and by stacking the real and imaginary parts of a complex-valued Jacobian on top of each other. The details about the employed FEM discretizations can be found in Section 5.

The necessary condition for minimizing (4.1) reads

$$(J^{(l)})^{\mathrm{T}}\Gamma^{-1}J^{(l)}\beta + a\begin{bmatrix}(\nabla R)(\varsigma)\\0\end{bmatrix} + b\begin{bmatrix}0\\(\nabla R)(\upsilon)\end{bmatrix} = (J^{(l)})^{\mathrm{T}}\Gamma^{-1}y^{(l)}, \qquad (4.2)$$

where $0 \in \mathbb{R}^N$. The gradient of $R : \mathbb{R}^N \to \mathbb{R}_+$ has the representation [15]

$$(\nabla R)(u) = H(u)u,$$

where $u = \varsigma$ or $u = v$ and

$$H_{m,l}(u) := \int_\Omega \frac{r'(|\nabla u(x)|)}{|\nabla u(x)|} \nabla \phi_m(x) \cdot \nabla \phi_l(x)\, dx, \qquad m, l = 1, \ldots, N. \tag{4.3}$$

This matrix can be interpreted as the FEM system matrix for the elliptic partial differential operator

$$-\nabla \cdot c_u \nabla \tag{4.4}$$

with the positive-valued diffusion coefficient

$$c_u : x \mapsto \frac{r'(|\nabla u(x)|)}{|\nabla u(x)|}, \quad \Omega \to \mathbb{R}_+.$$

Since the nodes on $S \subset \partial\Omega$ were excluded in (3.2), the stiffness matrix $H(u)$ corresponds to a natural boundary condition on $\partial\Omega \setminus \overline{S}$ and a homogeneous Dirichlet condition on $S$ for (4.4). In particular, both $H(\varsigma)$ and $H(v) \in \mathbb{R}^{N \times N}$ are positive definite under the reasonable assumption that there is a positive number FEM nodes on $S$.

The equation (4.2) can now be rewritten in the form

$$\left( (J^{(l)})^{\mathrm{T}} \Gamma^{-1} J^{(l)} + \begin{bmatrix} a\, H(\varsigma) & 0 \\ 0 & b\, H(v) \end{bmatrix} \right) \beta = (J^{(l)})^{\mathrm{T}} \Gamma^{-1} y^{(l)}, \tag{4.5}$$

which remains nonlinear since the matrices $H(\varsigma)$ and $H(v)$ depend on $\beta = [\varsigma^{\mathrm{T}}, v^{\mathrm{T}}]^{\mathrm{T}}$. We deal with this nonlinearity by, instead of solving (4.5) exactly, settling with taking a single step of lagged diffusivity iteration as in [15]: The $H$-matrices appearing in (4.5) are evaluated at the current iterate $\beta^{(l)} = [(\varsigma^{(l)})^{\mathrm{T}}, (v^{(l)})^{\mathrm{T}}]^{\mathrm{T}}$ and the resulting linear equation is (approximately) solved for the next one $\beta^{(l+1)}$. More precisely, denoting a Cholesky factor of $\Gamma^{-1}$ by $\Gamma^{-1/2}$ and setting

$$A = \Gamma^{-1/2} J^{(l)}, \qquad H = \begin{bmatrix} H(\varsigma^{(l)}) & 0 \\ 0 & \frac{b}{a} H(v^{(l)}) \end{bmatrix}, \qquad \tilde{y} = \Gamma^{-1/2} y^{(l)}, \tag{4.6}$$

we obtain the normal equation

$$(A^{\mathrm{T}} A + a H)\beta = A^{\mathrm{T}} \tilde{y}, \qquad a > 0, \tag{4.7}$$

to which we apply LSQR-aided priorconditioning accompanied by the Morozov discrepancy principle, as explained in the following subsection.

**4.2. Priorconditioning.** Since $H \in \mathbb{R}^{2N \times 2N}$ is positive definite, the unique solution of (4.7) is also the minimizer of the standard, quadratic Tikhonov functional

$$|A\beta - \tilde{y}|^2 + a\, \beta^{\mathrm{T}} H \beta.$$

From the Bayesian standpoint, this can be interpreted as having small $H$-norm of the parameter vector $\beta$ as prior information: Assigning to $H$ the role of a (scaled) inverse covariance matrix of a zero-mean Gaussian prior and assuming a suitable additive Gaussian measurement noise model, the solution of (4.7) is the *conditional mean* or *maximum a posteriori* estimate for the the linear model

$$A\beta = \tilde{y}. \tag{4.8}$$

We aim at preconditioning (4.7) so that the prior information in $H$ is included directly in the Krylov subspace structure produced by our preferred iterative solver, LSQR. The 'regularization parameter' $a > 0$ becomes superfluous in the process; however, the ratio $b/a$ of the free parameters appearing in the priors (3.5) is included in the definition of $H$ and must thus be chosen by the operator of the algorithm.

Since $H$ is positive definite, it allows a Cholesky decomposition $H = L^{\mathrm{T}} L$; in fact, any symmetric factorization of $H$ is adequate for our purposes and, even more importantly, it need not be formed explicitly in the reconstruction algorithm as we resort to the modified LSQR algorithm introduced in [1]. However, one needs to be able to apply $H^{-1}$ on a given vector. Multiplying (4.7) from the left by $(L^{-1})^{\mathrm{T}}$ and setting $a$ to zero yields

$$(L^{-1})^{\mathrm{T}} A^{\mathrm{T}} A L^{-1} \tilde{\beta} = (L^{-1})^{\mathrm{T}} A^{\mathrm{T}} \tilde{y} \qquad (4.9)$$

where $\tilde{\beta} = L\beta$. Although choosing $a = 0$ obviously makes the matrix equation (4.9) severely illconditioned, the prior information in $H$ is still included via the symmetric preconditioning by $L$.

We solve (4.9) by combining LSQR iteration [1] with an early stopping rule based on the Morozov discrepancy principle. Starting the iteration from $\tilde{\beta} = 0$, this means that after $m \in \mathbb{N}$ iterations the estimate for the logarithmic parameter vector $\beta$ belongs to the Krylov subspace

$$\mathcal{K}_m = \mathrm{span}\{H^{-1} A^{\mathrm{T}} \tilde{y}, (H^{-1} A^{\mathrm{T}} A) H^{-1} A^{\mathrm{T}} \tilde{y}, \ldots, (H^{-1} A^{\mathrm{T}} A)^{m-1} H^{-1} A^{\mathrm{T}} \tilde{y}\}$$

which clearly takes into account the prior information in $H$, that is, the approximate solution is in the range of $H^{-1}$ regardless of the number of iterations $m$. The LSQR iteration is terminated and the corresponding iterate dubbed $\beta^{(l+1)}$ when the residual corresponding to (4.8) reaches the (whitened) noise level

$$\varepsilon = \sqrt{\mathbb{E}\left(|\Gamma^{-1/2}\eta|^2\right)} = \sqrt{2KJ}$$

scaled by a suitable 'fudge factor' to account for unavoidable numerical noise.

**4.3. Implementation.** Collecting the material in the preceding two subsections, we have altogether introduced the following algorithm for simultaneous reconstruction of $\kappa$ and $\mu$.

ALGORITHM 1. *Select $T > 0$ in (3.7), the ratio $b/a$ for (3.5), a 'fudge factor' $\tau \geq 1$, and $S \subset \partial\Omega$ to support a homogeneous Dirichlet boundary condition for $\varsigma$ and $\upsilon$ (cf. (3.2), (4.3) and (4.4)). Let $\mathbf{1} = [1,\ldots,1]^{\mathrm{T}} \in \mathbb{R}^N$ and determine $(\kappa_0, \mu_0)$ as the minimizing pair for*

$$\left|\Gamma^{-1/2}\bigl(\mathcal{V} - \mathcal{M}(\log(\kappa)\mathbf{1}, \log(\mu)\mathbf{1})\bigr)\right|$$

*over $(\kappa, \mu) \in \mathbb{R}^2_+$. Set $l = 0$, $\epsilon = \sqrt{2KJ}$, $\beta^{(0)} = 0 \in \mathbb{R}^{2N}$ and $\mathcal{M} = \mathcal{M}(\beta^{(0)})$.*

1. *Evaluate the Jacobian of the map $\beta \mapsto \mathcal{M}(\beta)$ at $\beta^{(l)}$ and denote it by $J$.*

2. *Set $A = \Gamma^{-1/2} J$ and $\tilde{y} = \Gamma^{-1/2}\bigl(\mathcal{V} - \mathcal{M} + J\beta^{(l)}\bigr)$.*

3. *Build $H(\varsigma^{(l)})$ and $H(\upsilon^{(l)})$ according to (4.3) and form $H = L^{\mathrm{T}} L$ as in (4.6). (Recall that $\beta^{(l)} = [(\varsigma^{(l)})^{\mathrm{T}}, (\upsilon^{(l)})^{\mathrm{T}}]^{\mathrm{T}}$ and the factor $L$ is introduced only for notational convenience.)*

   *4. Apply the LSQR algorithm of [1] to*

$$(L^{-1})^{\mathrm{T}} A^{\mathrm{T}} A L^{-1} \tilde{\beta} = (L^{-1})^{\mathrm{T}} A^{\mathrm{T}} \tilde{y}, \qquad \beta = L^{-1} \tilde{\beta},$$

   *starting from $\tilde{\beta} = 0$. Denoting the LSQR sequence as $\{\beta_m^{(l)}\}_{m \in \mathbb{N}_0}$, terminate the iteration as soon as $\|A\beta_m^{(l)} - \tilde{y}\| \leq \tau \epsilon$. Denote the corresponding solution as $\beta^{(l+1)}$.*

   *5. Compute $\mathcal{M} = \mathcal{M}(\beta^{(l+1)})$. If*

$$\left|\Gamma^{-1/2}\bigl(\mathcal{V} - \mathcal{M}(\beta^{(l+1)})\bigr)\right| \leq \tau \epsilon$$

   *substitute $\beta^{(l+1)} = [(\varsigma^{(l+1)})^{\mathrm{T}}, (v^{(l+1)})^{\mathrm{T}}]^{\mathrm{T}}$ in (3.1) and declare the resulting $\kappa$ and $\mu$ the reconstruction. Otherwise, set $l \leftarrow l + 1$ and return to step 1.*

The choice of the free parameter $\tau$ in Algorithm 1 is case-dependent, but throughout our numerical experiments we use $S = \{\cup_k s_k\} \cup \{\cup_j m_j\}$ and the values $T = 5 \cdot 10^{-3}$ and $b/a = 1/3$ are fixed. The algorithm is insensitive to the choice of $S$ and $T > 0$, but selecting $b/a$ is a more delicate issue: According to our experience, sharp boundaries more easily appear in reconstructions of $\kappa$ compared to those of $\mu$, and we compensate for this phenomenon by assuming a looser prior, or less regularization for $\mu$, that is, we set $b/a < 1$ (cf. (3.5)).

**5. Numerical experiments.** To demonstrate the performance of Algorithm 1, we present four three-dimensional numerical experiments based on simulated data. In Case 1 we consider the situation, where we have only one unknown parameter to reconstruct (in turns, either $\kappa$ or $\mu$) as opposed to Case 2, where we reconstruct the two unknowns simultaneously. In these first two experiments we use unmodulated data, i.e. set $\omega = 0$ in (2.2). As a result, the imaginary parts of the measurements $M_{jk}$ remain zero, and hence the length of the (real-valued) measurement vector $\mathcal{M}$ effectively decreases to $JK$ (cf. (3.3)), making the computations less expensive. For comparison, Case 3 presents the corresponding results for frequency modulated data. The first three test cases consider a cylindrical body with the sources and sensors attached to its side. To conclude, we evaluate the functionality of the algorithm in a ball-shaped domain with nonconvex target inclusions in Case 4.

In all tests the data is generated, as explained in the beginning of Section 4.1, with FEM using piecewise linear basis functions on a dense finite element mesh (102562 nodes and 558244 tetrahedrons in Cases 1–3, 116520 nodes and 607553 tetrahedrons in Case 4). Subsequently, the simulated exact measurements are corrupted by 1% of normally distributed additive noise: we choose $\eta_i \sim \mathcal{N}(0, \sigma_i^2)$ and $\sigma_i = 0.01|\mathcal{M}_i|$ for all $i = 1, \ldots, 2JK$ in (3.4), which, in particular, makes the noise covariance matrix $\Gamma \in \mathbb{R}^{2JK \times 2JK}$ diagonal. In the reconstruction algorithm, we ignore the measurements on the sensors closest to the active photon source. The reason for this is two-fold: according to our experience, the numerical noise is the highest on these sensors, and the diffusion approximation is arguably the least accurate close to the photon source; see, e.g., [21] and the references therein.

**Case 1: Separate reconstructions from unmodulated data.** The domain in Cases 1–3 is a cylinder $\Omega = D(0,1) \times (0,1) \subset \mathbb{R}^3$, where $D(x,r) \subset \mathbb{R}^2$ denotes an open disk of radius $r > 0$ centered at $x \in \mathbb{R}^2$. There are $K = 24$ circular photon sources and $J = 24$ measurement sensors of the same size and shape attached evenly and alternatingly in three rings to the side of $\Omega$ (see the top row in Figure 5.1).

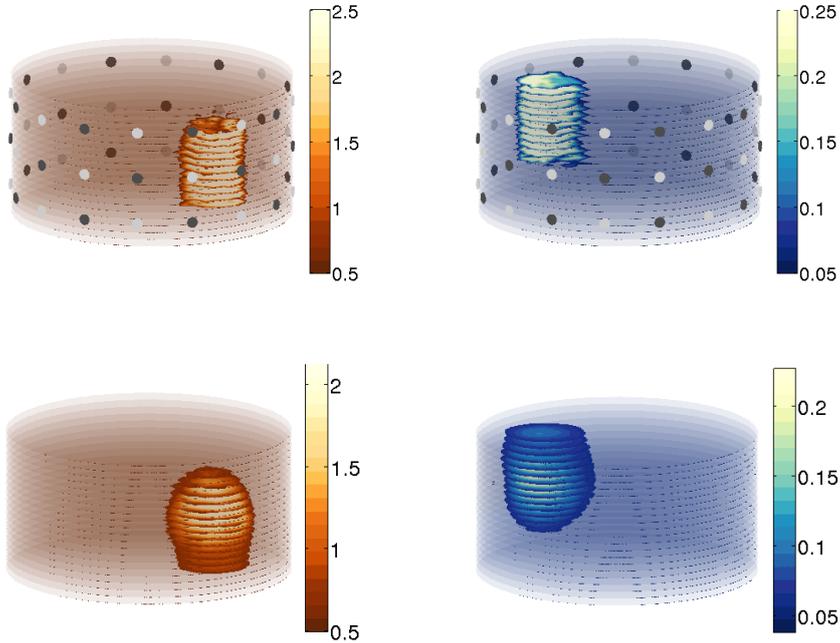

Fig. 5.1: Case 1. The separate reconstructions of absorption and diffusivity from unmodulated data when the value of the other parameter is a known constant. Top row: the target absorption (left) and diffusivity (right). The photon sources are plotted in light gray color and the measurement sensors in dark gray. Bottom row: the reconstructed absorption (left) and diffusivity (right). The values between $\mu_0 \pm 0.1$ and $\kappa_0 \pm 0.01$, respectively, are transparent.

Corresponding to each source, we simulate the measurements on all sensors, expect for the ones closest to the source in question, which amounts to three or four passive sensors depending on the position of the active source. For unmodulated data, this leads to a real-valued measurement vector of length 496.

We begin by considering the situation, where either the diffusivity or the absorption is homogeneous with a known constant value, and hence we only need to reconstruct the other material parameter. We first fix the diffusivity to the (known) constant value $\kappa(x) = 0.05$ for all $x \in \Omega$. The corresponding target absorption is illustrated in the top left image of Figure 5.1; the absorption equals 0.5 apart from a cylindrical inclusion of radius 0.2 and height 0.6 touching the bottom of $\Omega$ and having a higher absorption level $\mu = 2.5$. In the second part of this first test case, the absorption is in turn fixed to the (known) constant value $\mu(x) = 0.5$ for all $x \in \Omega$. The associated target diffusivity is shown in the top right image of Figure 5.1; the diffusion coefficient takes the constant value $\kappa = 0.05$ except for a cylindrical inclusion of radius 0.2 and height 0.6 touching the top of $\Omega$ and exhibiting a higher level of diffusivity $\kappa = 0.25$. (The background parameter levels of $\mu = 0.5$ and $\kappa = 0.05$ in the 'unit cylinder' $\Omega$ could model, e.g., a real-world cylinder of radius 1 cm and height 1 cm with constant absorption $0.5 \, \text{cm}^{-1}$ and reduced scattering coefficient $6.17 \, \text{cm}^{-1}$, cf. [2, 27].)

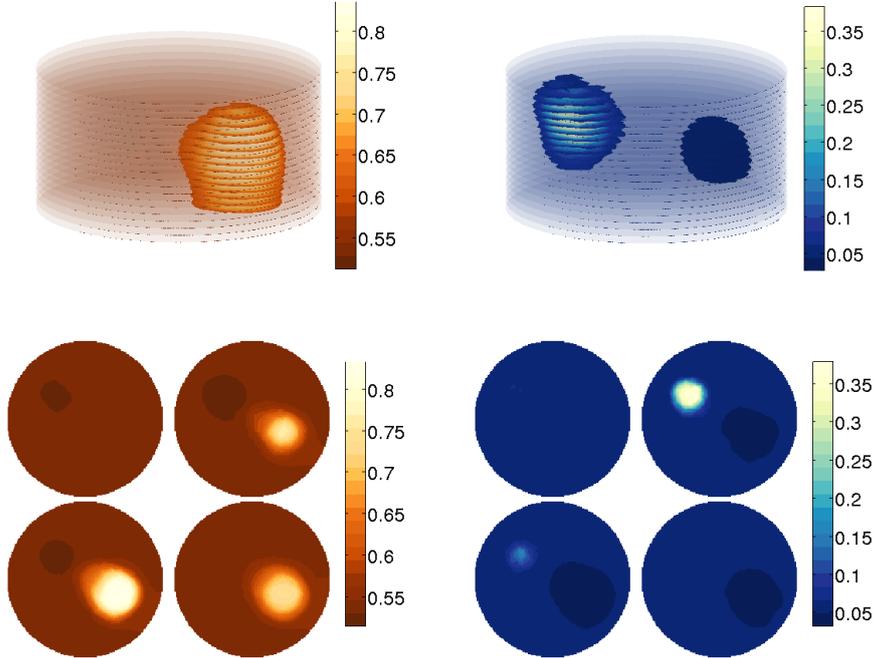

Fig. 5.2: Case 2. The simultaneous reconstruction of absorption and diffusivity from unmodulated data. Top row: the reconstructed absorption (left) and diffusivity (right). The values between $\mu_0 \pm 0.1$ and $\kappa_0 \pm 0.01$, respectively, are transparent. Bottom row: horizontal slices of the reconstruction at heights 0.9, 0.6, 0.4 and 0.1.

The separate reconstructions of the absorption and the diffusivity are presented in the bottom left and bottom right images of Figure 5.1, respectively. To enhance the visualizations, the values close to the reconstructed background parameter levels $\mu_0 = 0.55$ and $\kappa_0 = 0.052$ are transparent in the respective images; see the initialization phase of Algorithm 1. Neither of the reconstructed inclusions is quite of the correct shape; the algorithm (or the prior) clearly prefers oval objects over cylindrical ones. Furthermore, the parameter levels within the reconstructed inclusions are slightly too low, though the maximal reconstructed values almost reach the target levels: $\mu_{\max} = 2.12$ and $\kappa_{\max} = 0.227$. In both reconstructions, the respective inclusion appears roughly in the correct position, it is well localized and the background is approximately constant, which are the features the algorithm was designed to promote. All in all, the separate reconstructions are comparable to those presented in [15] for EIT.

To compensate for numerical errors, the fudge factor was set to $\tau = 1.3$ in both subcases. The homogeneous estimates for the background parameter levels were computed on a coarse finite element mesh with 32 130 nodes and 161 539 tetrahedrons. Subsequently, the actual reconstructions were formed on a mesh with 56 021 nodes and 296 276 tetrahedrons. The running time of Algorithm 1 was approximately 5.5 minutes for reconstructing the absorption (convergence after four linearizations) and about 3 minutes for the diffusivity (convergence after two linearizations) with a MATLAB (2014a) implementation run on a laptop with 16 GB RAM and an Intel Core i7-4600U CPU having clock speed 2.10 GHz.

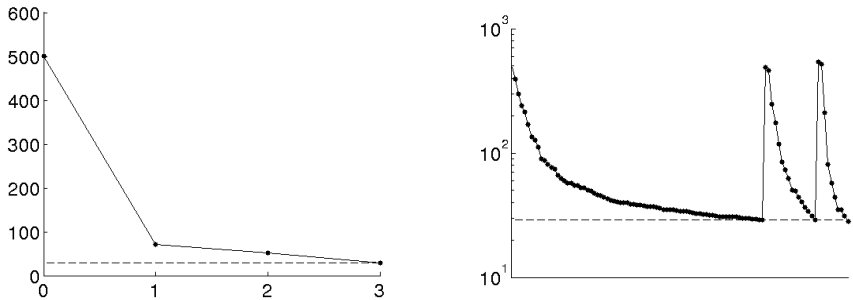

Fig. 5.3: Case 2. The convergence of the algorithm for the simultaneous reconstruction of absorption and diffusivity from unmodulated data. Left: The residuals $|\Gamma^{-1/2}(\mathcal{V} - \mathcal{M}(\beta^{(l)}))|$ after each linearization. Right: The residuals $|A\beta_m^{(l)} - \tilde{y}|$ after each LSQR step (on a logarithmic scale). The dashed lines indicate the target residual level $\tau\epsilon$.

**Case 2: Simultaneous reconstruction from unmodulated data.** In the second experiment, we combine the two separate parts of Case 1 and aim at reconstructing the absorption and diffusivity simultaneously. The domain and the target coefficients are the same as in Case 1, that is, the to-be-reconstructed absorption and diffusivity are as shown in the top row of Figure 5.1. We still consider unmodulated measurements, i.e. set $\omega = 0$.

The results produced by Algorithm 1 are presented in Figure 5.2. The top row illustrates the three-dimensional reconstructions with transparency applied to parameter values close to the estimated background levels $\mu_0 = 0.51$ and $\kappa_0 = 0.049$, and the bottom row shows horizontal slices of the reconstructions. As in Case 1, the inclusions do not have correct shapes, but the background is almost constant and the localization of the inhomogeneities is good. However, when reconstructing the absorption and diffusivity simultaneously, the parameter levels within the inclusions are not as accurate as in Case 1, and we also observe some 'cross-talk' between the two parameters. More precisely, the algorithm seems to partially explain the increase in the absorption by decreasing the diffusivity, and there is also a small drop in the absorption at the location of the target inclusion with high diffusivity. The ratio parameter $b/a$ controls to a certain extent which of the two unknowns the algorithm is more prone to change in order to explain the measurements. Nevertheless, tuning this ratio does not generally remove the observed cross-talk effect.

Figure 5.3 illustrates the convergence of the algorithm. The left-hand image depicts the residuals after each linearization of the forward model, which corresponds to the overall stopping criterion in step 5 of Algorithm 1. The right-hand image shows the residual for the linear model (4.8) after each LSQR step, which in turn monitors the termination of the inner loop in step 4 of Algorithm 1. As indicated by Figure 5.3, the algorithm converges after three rounds of LSQR iterations, requiring fewer LSQR steps each round. We conclude that the algorithm functions stably and produces results that are compatible with the priors. This suggests that the cross-talk between the absorption and the diffusivity in the reconstructions of Figure 5.2 may be an unavoidable consequence of the fundamental nonuniqueness in diffuse optical tomography [4]. This phenomenon was also encountered, e.g., in [29].

The computations in Case 2 were performed on the same finite element meshes and with the same hardware as in Case 1. We also used the same fudge factor $\tau = 1.3$. The running time of the algorithm was approximately 5 minutes.

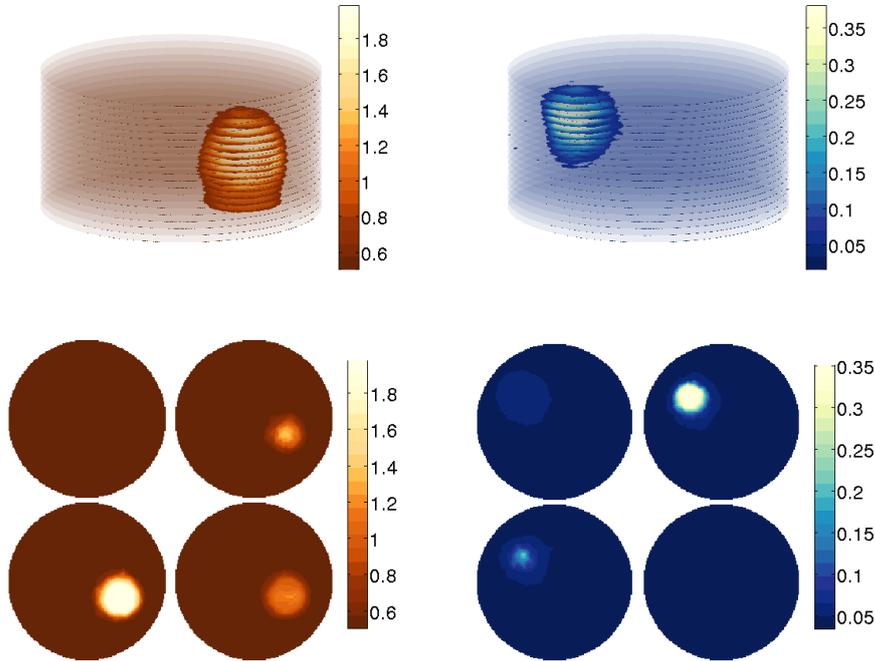

Fig. 5.4: Case 3. The simultaneous reconstruction of absorption and diffusivity from frequency modulated data. Top row: the reconstructed absorption (left) and diffusivity (right). The values between $\mu_0 \pm 0.1$ and $\kappa_0 \pm 0.01$, respectively, are transparent. Bottom row: horizontal slices of the reconstruction at heights 0.9, 0.6, 0.4 and 0.1.

**Case 3: Simultaneous reconstruction from frequency modulated data.** In order to reduce the cross-talk between the absorption and diffusivity, we next resort to frequency modulated data. We still consider the domain and target material parameters shown in the top row of Figure 5.1, but now the inward photon flux at the boundary sources is harmonically modulated so that $\omega/c \approx 0.021$. This corresponds approximately to the frequency of 100 MHz if the unit of length in $\Omega$ is cm. Since the imaginary parts of the measurements $M_{jk}$ are now nonzero, the length of the real-valued measurement vector $\mathcal{M}$ increases to 992 (cf. (3.3)).

The reconstruction is visualized in Figure 5.4. Once again, we observe that Algorithm 1 produces results that are consistent with the prior: the reconstructions of the absorption and diffusivity in Figure 5.4 are composed of a homogeneous background with well localized, yet somewhat incorrectly shaped, inclusions at approximately correct locations. Even more importantly, the employment of frequency modulated data indeed helped us to get rid of the unwanted cross-talk between the absorption and diffusivity. The parameter levels in the inclusions are still slightly off, but they mimic the target considerably better than in Case 2. The convergence of the algorithm is reported in Figure 5.5 which is organized in the same way as Figure 5.3. This time the algorithm converges after four LSQR rounds and, as in Case 2, the number of the needed LSQR steps decreases round by round.

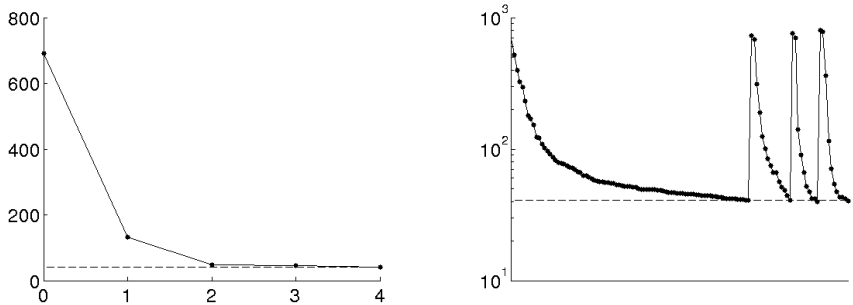

Fig. 5.5: Case 3. The convergence of the algorithm for the simultaneous reconstruction of absorption and diffusivity from frequency modulated data. Left: The residuals $|\Gamma^{-1/2}(\mathcal{V} - \mathcal{M}(\beta^{(l)}))|$ after each linearization. Right: The residuals $|A\beta_m^{(l)} - \tilde{y}|$ after each LSQR step (on a logarithmic scale). The dashed lines indicate the target residual level $\tau\epsilon$.

The employed finite element meshes were again the same as previously and, as before, we used the fudge factor $\tau = 1.3$. In this case the measurement vector $\mathcal{M}$ was twice as long as in Cases 1 and 2, which slowed down the algorithm. The running time was approximately 11 minutes with the same hardware as before. The homogeneous approximation for the absorption was $\mu_0 = 0.55$ and for the diffusion coefficient $\kappa_0 = 0.051$.

**Case 4: Ball-shaped domain.** In the final experiment, the domain $\Omega$ is a ball of radius 10 with $K = 32$ circular photon sources and $J = 60$ circular measurement sensors distributed evenly on its surface (see the top row of Figure 5.6). For each source, we simulate the measurements on all sensors apart from the five or six ones lying the closest to the source in question. This results in a real-valued measurement vector of length 3480.

The target material parameters are illustrated in the top row of Figure 5.6. The absorption level of $\Omega$ is $\mu \equiv 0.025$ except for an inclusion where the absorption increases to $\mu = 0.125$. On the other hand, the diffusivity satisfies $\kappa \equiv 0.15$ apart from an inhomogeneity where, in contrast to Cases 1–3, it decreases to $\kappa = 0.075$. (The background levels of $\mu = 0.025$ and $\kappa = 0.15$ in the ball $\Omega$ could model, e.g., a real-world ball of radius $1\,\text{cm}$ with constant absorption $0.25\,\text{cm}^{-1}$ and reduced scattering coefficient $22.0\,\text{cm}^{-1}$, cf. [2, 27].) Both aforementioned inclusions are nonconvex: they have three-dimensionally twisted z-shapes such that the three sections are parallel to $x$, $y$ and $z$-axis, respectively. In particular, the middle sections of the inclusions (lying at height $z = 2$) coincide, meaning that there is a jump in both material parameters at the same location. The measurements are modulated so that $\omega/c \approx 0.0126$, which corresponds approximately to the frequency of $600\,\text{MHz}$ if the unit of length in $\Omega$ is mm.

A three-dimensional visualization of the reconstruction is presented in the middle row of Figure 5.6, while the bottom row shows corresponding horizontal cross sections. The properties of the reconstruction are similar to the previous case: the positions and localization of the inclusions are good, the background is approximately constant and, in particular, the nonconvexity of the inclusions is captured reasonably well. In addition, even though the inclusions coincide in their middle sections, the cross-talk between them is minimal. However, it should be noted that it is necessary to use a fairly high angular frequency $\omega$ to achieve this. As before, the parameter values within the reconstructed inclusions do not quite reach the target levels. Some instability is also observed near the boundary of $\Omega$, which causes small jumps in the reconstructed diffusivity. Although these jumps occur close to the photon sources and the measurement sensors, according to our experiments, neither changing the choice of $S \subset \partial\Omega$ (cf. (3.2)) nor using a natural boundary condition for (4.4) everywhere on $\partial\Omega$ removes the artefacts.

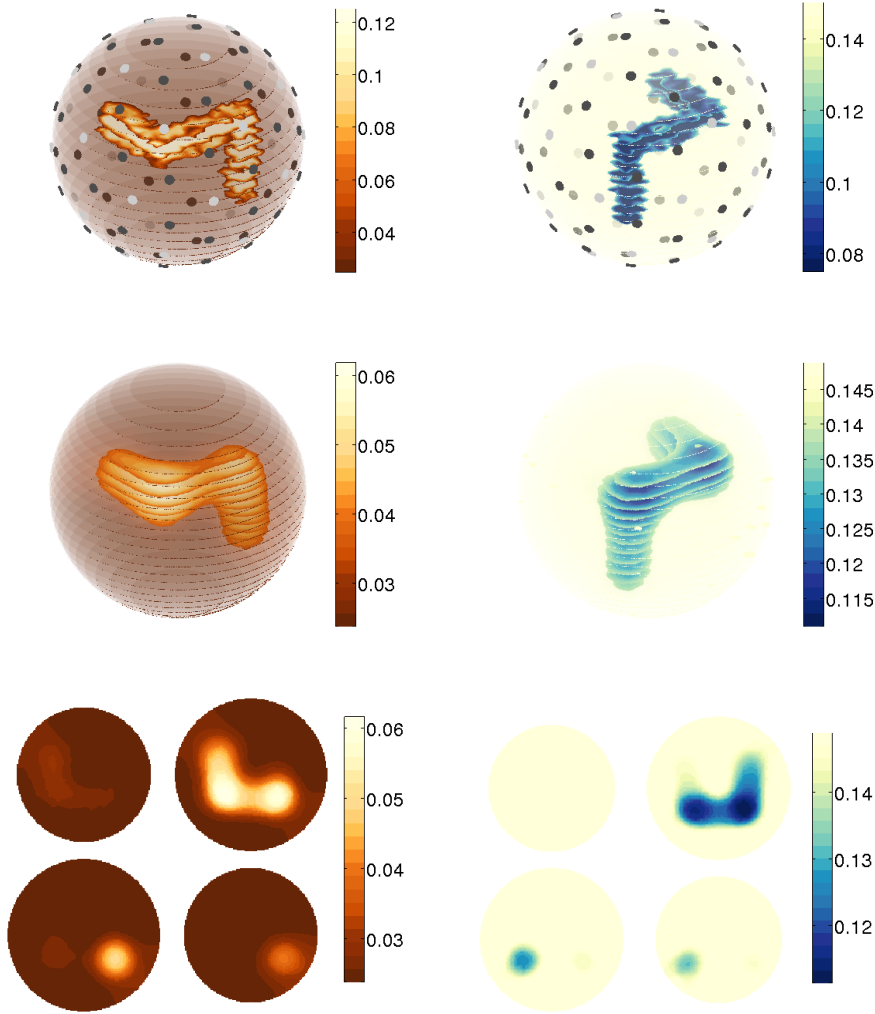

Fig. 5.6: Case 4. The simultaneous reconstruction of absorption and diffusivity from frequency modulated data. Top row: the target absorption (left) and diffusivity (right). The photon sources are plotted in light gray color and the measurement sensors in dark gray. Middle row: the reconstructed absorption (left) and diffusivity (right). The values between $\mu_0 \pm 0.01$ and $\kappa_0 \pm 0.01$, respectively, are transparent. Bottom row: horizontal slices of the reconstructions at heights 5, 2, $-2$ and $-5$.

To ensure convergence, we used a somewhat higher fudge factor $\tau = 1.6$ in this final test case. The homogeneous estimates for the parameters $\mu_0 = 0.0282$ and $\kappa_0 = 0.149$ were computed on a coarse finite element mesh with 23648 nodes and 101305 tetrahedrons, whereas the actual reconstruction was formed on a mesh with 49072 nodes and 232750 tetrahedrons. The hardware was as described in Case 1. Even though both of the aforementioned meshes were slightly coarser than the corresponding ones in Cases 1–3, the more than three times higher number of measurements resulted in a slightly longer run time of 13.5 minutes compared to Case 3. The algorithm converged after four linearizations of the measurement model.

**6. Concluding remarks.** We have introduced an edge-enhancing technique for reconstructing the absorption and diffusion coefficients in diffuse OT. The method is a variant of the algorithm introduced in [15] for EIT; see also [1]. By means of three-dimensional numerical experiments with simulated data, we have demonstrated that the algorithm is capable of simultaneously locating inclusions in both the absorption and the diffusivity, if the inward photon flux through the boundary sources is modulated with a high enough harmonic frequency and if the to-be-reconstructed inhomogeneities in the material parameters exhibit sufficient contrasts compared to the background. The running time of the algorithm on a standard laptop computer is around ten minutes for 1000 photon flux measurement at the boundary and a parametrization of the unknowns with $10^5$ degrees of freedom.

**Acknowledgments.** The work of A. Hannukainen, N. Hyvönen and H. Majander was supported by the Academy of Finland (decisions 267297 and 267789).

REFERENCES

[1] Arridge, S., Betcke, M., and Harhanen, L. Iterated preconditioned LSQR method for inverse problems on unstructured grids. *Inverse Problems 30* (2014), 075009.
[2] Arridge, S. R. Optical tomography in medical imaging. *Inverse Problems 15* (1999), R41–R93.
[3] Arridge, S. R., and Hebden, J. C. Optical imaging in medicine: II. Modelling and reconstruction. *Phys. Med. Biol. 42* (1997), 825–840.
[4] Arridge, S. R., and Lionheart, W. R. B. Nonuniqueness in diffusion-based optical tomography. *Opt. Lett. 23* (1998), 882–884.
[5] Boas, D. A., Brooks, D. H., Miller, E. L., DiMarzio, C. A., Kilmer, M., Gaudette, R. J., and Zhang, Q. Imaging the body with diffuse optical tomogrphy. *IEEE Signal Process. Mag.* (2001), 57–75.
[6] Calvetti, D. Preconditioned iterative methods for linear discrete ill-posed problems from a bayesian inversion perspective. *J. Comput. Appl. Math. 198* (2007), 378395.
[7] Calvetti, D., McGivney, D., and Somersalo, E. Left and right preconditioning for electrical impedance tomography with structural information. *Inverse Problems 28* (2012), 055015.
[8] Calvetti, D., and Somersalo, E. Priorconditioners for linear systems. *Inverse problems 21* (2005), 1397–1418.
[9] Calvetti, D., and Somersalo, E. Hypermodels in the bayesian imaging framework. *Inverse Problems 24* (2008), 034013.
[10] Dierkes, T., Dorn, O., Natterer, F., Palamodov, V., and Sielschott, H. Fréchet derivatives for some bilinear inverse problems. *SIAM J. Appl. Math. 62* (2002), 2092–2113.
[11] Eldén, L. A weighted pseudoinverse, generalized singular values, and constrained least squares problem. *BIT 22* (1982), 487–502.
[12] Gibson, A. P., Hebden, J. C., and Arridge, S. R. Recent advances in diffuse optical imaging. *Phys. Med. Biol. 50* (2005), R1–R43.
[13] Grisvard, P. *Elliptic Problems in Nonsmooth Domains*. Pitman, 1985.
[14] Hansen, P. C. *Rank-Deficient and Discrete Ill-Posed Problems*. SIAM, 1998.
[15] Harhanen, L., Hyvönen, N., Majander, H., and Staboulis, S. Edge-enhancing reconstruction algorithm for three-dimensional electrical impedance tomography. *SIAM J. Sci. Comput. 37* (2015), B60–B78.

[16] HARRACH, B. On uniqueness in diffuse optical tomography. *Inverse Problems 25* (2009), 055010.
[17] HEBDEN, J. C., ARRIDGE, S. R., AND DELPY, D. T. Optical imaging in medicine: I. Experimental techniques. *Phys. Med. Biol. 42* (1997), 825–840.
[18] HEINO, J., AND SOMERSALO, E. Estimation of optical absorption in anisotropic background. *Inverse Problems 18* (2002), 559–573.
[19] HILGERS, J. W. On the equivalence of regularization and certain reproducing kernel hilbert space approaches for solving first kind problems. *SIAM J. Numer. Anal. 13* (1976), 172–184.
[20] HINTERMÜLLER, M., AND WU, T. Nonconvex TV$^q$ -models in image restoration: Analysis and a trust-region regularizationbased superlinearly convergent solver. *SIAM J. Imag. Sci. 6* (2013), 13851415.
[21] LEHTIKANGAS, O., TARVAINEN, T., AND KIM, A. D. Modeling boundary measurements of scattered light using the corrected diffusion approximation. *Biomed. Opt. Express 3* (2012), 552–571.
[22] NISSILÄ, I., NOPONEN, T., KOTILAHTI, K., KATILA, T., LIPIÄINEN, L., TARVAINEN, T., SCHWEIGER, M., AND ARRIDGE, S. Instrumentation and calibration methods for the multichannel measurement of phase and amplitude in optical tomography. *Rev. Sci. Instrum. 76* (2005), 044302.
[23] PAIGE, C. C., AND SAUNDERS, M. A. Algorithm 583: LSQR: Sparse linear equations and least squares problems. *ACM Trans. Math. Softw. 8* (1982), 195–209.
[24] PAIGE, C. C., AND SAUNDERS, M. A. LSQR: An algorithm for sparse linear equations and sparse least squares. *ACM Trans. Math. Softw. 8* (1982), 43–71.
[25] PERONA, P., AND MALIK, J. Scale-space and edge detection using anisotropic diffusion. *IEEE T. Pattern Anal. 12* (1990), 629639.
[26] RUDIN, L. I., OSHER, S., AND FATEMI, E. Nonlinear total variation based noise removal algorithms. *Physica D 60* (1992), 259–268.
[27] SANDELL, J. L., AND ZHU, T. C. A review of in-vivo optical properties of human tissues and its impact on pdt. *J. Biophotonics 4* (2011), 773–787.
[28] SCHWEIGER, M., ARRIDGE, S. R., AND NISSILÄ, I. Gauss-newton method for image reconstruction in diffuse optical tomography. *Phys. Med. Biol. 50* (2005), 2365–2386.
[29] SCHWEIGER, M., AND R., A. S. Application of temporal filters to time resolved data in optical tomography. *Phys. Med. Biol. 44* (1999), 1699–1717.
[30] TARVAINEN, T., VAUHKONEN, M., AND ARRIDGE, S. R. Gauss-newton reconstruction method for optical tomography using the finite element solution of the radiative transfer equation. *J. Quant. Spectrosc. Radiat. Transfer 109* (2008), 2767–2778.
[31] VOGEL, C. R., AND OMAN, M. E. Iterative methods for total variation denoising. *SIAM J. Sci. Comput. 17* (1996), 227238.
[32] WLOKA, J. *Partial Differential Equations*. Cambridge University Press, 1987.